\documentclass[14pt]{amsart}
\usepackage{amssymb}
\usepackage{amscd}
\usepackage[all]{xy}

\numberwithin{equation}{section}

\def\today{\number\day\space\ifcase\month\or   January\or February\or
   March\or April\or May\or June\or   July\or August\or September\or
   October\or November\or December\fi\   \number\year}

\theoremstyle{definition}

\newtheorem{thm}{Theorem}[section]
\newtheorem{lem}[thm]{Lemma}
\newtheorem{prp}[thm]{Proposition}
\newtheorem{dfn}[thm]{Definition}
\newtheorem{cor}[thm]{Corollary}

\newtheorem*{Thm*}{Theorem}

\newtheorem{claim}[thm]{Claim}

\newcommand{\beq}{\begin{equation}}
\newcommand{\eeq}{\end{equation}}
\newcommand{\beqr}{\begin{eqnarray*}}
\newcommand{\eeqr}{\end{eqnarray*}}
\newcommand{\bal}{\begin{align*}}
\newcommand{\eal}{\end{align*}}
\newcommand{\bei}{\begin{itemize}}
\newcommand{\eei}{\end{itemize}}

\newcommand{\af}{\alpha}

\newcommand{\ep}{\varepsilon}

\newcommand{\cF}{{\mathcal{F}}}

\newcommand{\cI}{{\mathcal{I}}}
\newcommand{\cJ}{{\mathcal{J}}}
\newcommand{\cK}{{\mathcal{K}}}

\newcommand{\cM}{{\mathcal{M}}}

\newcommand{\cS}{{\mathcal{S}}}

\newcommand{\cW}{{\mathcal{W}}}

 \def \xn{\{ x_n \}_{n=1}^{\infty}}

\def \wn{\{ w_n \}_{n=1}^{\infty}}

\def \hn{\{h_n \}_{n=1}^{\infty}}

\def \deltan{\{ \delta_n \}_{n=1}^{\infty}}

\def\pf{ \noindent {\bf Proof: \  }}

\def \e{\varepsilon}
\def \d{\delta}

\title{Approximate identities for Ideals in $L(L^p))$}

\author{W.~B.~Johnson$^*$\ and  \
G.~Schechtman}

\date{January 2024}

\address{}

\email{}


\thanks{$^*$Supported in part by NSF DMS-1900612}

\begin{document}

\maketitle

\begin{abstract}

The main result is that the only non trivial closed ideal in the Banach algebra $L(L^p)$ of
bounded linear operators on $L^p(0,1)$, $1\le p < \infty$, that has a left approximate identity is
the ideal of compact operators. The algebra $L(L^1)$ has at least one non trivial closed ideal
that has a contractive right approximate identity as well as many, including the unique maximal
ideal, that do not have a right approximate identity.
\end{abstract}

\section{Introduction}
Arguably the most natural  examples of $L^p$,
$1\le p < \infty$,  operator algebras; i.e.,
closed subalgebras of the Banach algebra $L(L^p)$ of (always bounded, linear) operators on
$L^p(\mu)$
for some measure $\mu$, are the closed ideals in $L(L^p)$. When $p=2$ and the measure $\mu$ is
separable,  it is well-known that the only non trivial closed ideal in $L(L^p)$ is the ideal
$\cK(L^p)$ of compact operators.
Similarly, $\cK(\ell^p)$ is the only non trivial closed ideal in $L(\ell^p)$ for other values of $1\le
p < \infty$. Here we focus on separable $L^p(\mu)$, in which case up to Banach space isomorphism the
only other examples are $L^p(0,1)$, so in the sequel we let $L^p$ denote $L^p(0,1)$.

When $1<p\not=2 <\infty$, the structure of the closed ideals in $L(L^p)$ is complicated and not well
understood. Given any Banach space $X$, every non trivial ideal in $L(X)$ contains the ideal
$\cF(X)$ of finite rank operators on $X$, and hence every non trivial closed ideal contains the
ideal $\cK(X)$ of compact operators on $X$ when $X$ has the approximation property; in particular,
when $X$ is any $L^p$ space, $1\le p \le \infty$.  Sometimes $L(X)$ has a largest non trivial ideal.
Let $\cM(X)$ denote the collection of all operators $T\in L(X)$ such that the identity operator
$I_X$ on $X$ does not factor through $T$ (an operator $S:W\to Z$ is said to factor through an
operator $T:X\to Y$ provided there exist operators $B : W\to X$ and $A: Y\to Z$ such that $S =
ATB$). It is easy to see that if $\cM(X)$ is closed under addition then $\cM(X)$ is an ideal, in
which case $\cM(X)$ is obviously the largest ideal in $L(L^p)$.  It is true but non trivial that
$\cM(L^p)$ is closed under addition for $1\le p < \infty$ (see \cite{es} for the case $p=1$ and
\cite{jmst} for the cases $1<p<\infty$).  Actually more is contained in the cited papers; namely,
that  $\cM(L^p)$ is the class of $L^p$--singular operators in $L(L^p)$. (An operator $T:X\to Y$ is
said to be {\sl $Z$--singular} provided there is no subspace $Z_1$ of $X$ such that $Z_1$ is
isomorphic to $Z$ and the restriction $T_{|Z_1}$ of $T$ to $Z_1$ is an isomorphism. So $T$ is {\sl
strictly singular} if it is $Z$-singular for every infinite dimensional $Z$.) This last description of 
$\cM(L^p)$ will be used later in this note.

We take this oportunity to describe some of the ways to build many more exmaples of closed ideals in $L(L^p)$.
As was pointed out by Pietsch \cite{pie} it follows from  \cite{sch}  that there are infinitely
many different  closed ideals in $L(L^p)$, $1<p\not= 2 <\infty$. In a similar way it follows from \cite{brs} that there 
are
at least $\aleph_1$. These ideals are all {\sl large} in the sense that they contain operators that
are not strictly singular.  Schlumprecht and  Zs\'ak \cite{sz} proved that there are at least a
continuum $c:=2^{\aleph_0}$ closed ideals in $L(L^p)$, $1<p\not= 2 <\infty$ . The ideals constructed
in \cite{sz} are all {\sl small}; i.e., contained in the strictly singular operators. Recently,
the authors \cite{js}
proved that there are $2^c$ large closed ideals and also $2^c$ small closed ideals in $L(L^p)$ when
$1<p\not=2 <\infty$.
The  ideals constructed in \cite{sch} and \cite{brs} look much different from the  ideals
constructed
 in \cite{sz} and \cite{js} and illustrate somewhat different approaches to constructing ideals.  In
 the first approach one   constructs a family
$\mathcal{S}$ of complemented subspaces of $L^p$ that are mutually non isomorphic (in fact, for each two members of 
$\mathcal{S}$ one of them does not embed complementably into any finite direct sum of the other). Given one of these 
spaces $X$, let $\mathcal{I}_X$ denote the space of all operators
in $L(L^p)$ that factor through some finite direct sum of $X$.
Then $\mathcal{I}_X$ is obviously a (possibly non closed) ideal in $L(L^p)$, and if an idempotent
$P$ is in the closure $\overline{\mathcal{I}_X}$ of $\mathcal{I}_X$, then it is easy to see that
the range of
$P$ is isomorphic to a complemented subspace of some finite direct sum of $X$.  Consequently,
$ \overline{\mathcal{I}_X} \not= \overline{\mathcal{I}_Y}$ for distinct $X$, $Y$ in $\mathcal{S}$.
In practice, it is a nuisance to have to consider all finite direct sums of the spaces $X$ and so
one builds examples where it is easy to check that $X\oplus X$ is isomorphic to $X$. In any case,
the ideals $ \overline{\mathcal{I}_X} $ are  large ideals because, by construction, they contain
projections onto infinite dimensional subspaces, and such operators are  obviously not strictly
singular.

In the second approach, one defines  a specific operator $T$ from some Banach space $X_1$ into some
Banach space $X_2$ and considers the collection  $\mathcal{I}_T$ of all operators $S$ in $L(L^p)$
that factor through $T$; i.e., $S $ can be written as $ASB$ for some bounded linear operators $B:
L^p \to X_1$ and $A: X_2 \to L^p$.  In practice one uses only $T$ for which the direct sum operator
$T\oplus T : X_1 \oplus X_1 \to X_2 \oplus X_2$ factors through $T$ in order to guarantee that
$\mathcal{I}_T$  is a (possibly non closed) ideal in $L(L^p)$ so that the closure
$\overline{\mathcal{I}_T}$ is a closed ideal in $L(L^p)$.  Of course, the first approach is the
special case of the second approach where $T$ is an idempotent in $L(L^p)$, but in the second
approach what have been used are  injections, including strictly singular operators such as  the
formal injections $i_{s,r}$
from $\ell_s$ to $\ell_r$ with $s<r$. When the injection $T$ is strictly singular, the deal
$\overline{\mathcal{I}_T}$ consists of strictly singular operators and so is  a  small ideal. 
In an earlier version of \cite{js} this approach was used to construct $c$ large closed ideals 
in $L(L^p)$, $1<p\not=2<\infty$. The
injections used there are not strictly singular but are very far
from being idempotents.

In \cite{jps}, using this approach but with very different operators, it was proved that  $L(L^1)$ 
has at least $c$ small closed ideals. It is open
whether  $L(L^1)$ has more than a continuum of closed ideals. There are known to exist only three
large ideals in $L(L^1)$.

To build $2^c$ closed ideals in $L(L^p)$, $1<p\not=2<\infty$, a yet somewhat different approach was 
used in \cite{js}. Each of these ideals is generated by a class of operators $\{T_\alpha\}_{\alpha\in \mathcal A}$. 
That is the closure of the set of operators of the form $\sum_{i=1}^n A_iT_{\alpha_i}B_i$ where n ranges 
over the positive integers and $\{\alpha_1,\dots,\alpha_n\}$ over ${\mathcal A}^n$. This of course is not 
very informative as any ideal is of that form. The main contribution of \cite{js} is to find a way to produce 
$2^c$ sets ${\mathcal A}$ that give rise to different ideals. 

Although the classification of the closed ideals in $L(L^p)$ is far from being known, we can prove
some general facts about them. These facts involve the notion of approximate identity.

\begin{dfn} Given a Banach space $X$ and an ideal $\cI$
in $L(X)$, we say that a net $(T_\af)$ of operators
in $\cI$ is a left; respectively,  right,  approximate identity for $\cI$ provided
 $\| T_{\af} S - S \| \to 0$; respectively,  $\| ST_{\af}  - S \| \to
0$,   for every $S$ in $\cI$.
\end{dfn}

In particular, we prove in Section \ref{left}  that the only non
trivial ideal in $L(L^p)$, $1\le p < \infty$, that has a left approximate identity is the ideal
$\cK(L^p)$ of compact operators.  If follows by duality that for $1<p<\infty$ the ideal $\cK(L^p)$ is
the  only non trivial ideal that has a right approximate identity. In contrast to this last fact, in
Section \ref{L^1} we prove that $L(L^1)$ has a non trivial ideal other than $\cK(L^1)$ that has a
right approximate identity consisting of contractive idempotents.  We also point out that no small
ideal in $L(L^1)$ has a right  approximate identity and  show that at least one of the other two
known large ideals in $L(L^1)$, the unique maximal ideal $\cM(L^1)$, does not have  a right  approximate identity.

We thank N.~C.~Phillips for many discussions about the contents herein. In fact, Phillips raised the
question whether there are  ideals in $L(L^p)$ other than $\cK(L^p)$ that have a left or right
approximate identity. Part of his motivation was the open problem to determine which ideals $\cI$ in
$L(L^p)$ have the property that $L(L^p)/\cI$ are isomorphically homomorphic to a subalgebra of
$L(L^p(\mu))$ for some measure $\mu$. The proofs in \cite{bj} and \cite{bp} that $\cK(L^p)$ has this
property both use the fact that $\cK(L^p)$ has a bounded approximate identity.

\section{Left approximate identities for ideals in $L(L^p)$}\label{left}

Here we prove

\begin{thm} \label{main}
Let $\mathcal{I}$ be a non trivial closed ideal in $L(L^p)$, $1\le p < \infty$, and assume that
$\mathcal{I}$ has a left approximate identity. Then $\mathcal{I}$ is $\cK(L^p)$, the ideal of compact
operators.
\end{thm}

\pf Since the Haar basis is a monotone Schauder basis for $L^p$,   the ideal $\cK(L^p)$ has a left
approximate identity consisting of finite rank projections of norm one. The converse direction
really consists of three theorems as the proofs for the cases $p=1$, $2<p<\infty$, and $1<p<2$  are
very different.  To call attention to the fact that we prove something stronger than an ideal has no
left approximate identity, we introduce some new language.  Given a Banach space $X$, an ideal $\cI$
in $L(X)$, and a collection $\cS$ of operators in $L(X)$, we say that a net $(T_\af)$ of operators
in $L(X)$ is an $\cS$--left; respectively, $\cS$--right,  approximate identity for $\cI$ provided
each $T_\af$ is in $\cS$ and $\| T_{\af} T - T \| \to 0$; respectively,  $\| TT_{\af}  - T \| \to
0$,   for every $T$ in $\cI$.

\bigskip\noindent
Case $1$: $p=1$.

\bigskip
Let $\mathcal{I}$ be a non trivial closed ideal in $L(L^1)$.  We treat first the case  that  $\cI$
is contained in the closed ideal $\mathcal{W} =
\mathcal{W}(L^1)$ of weakly compact operators in $L(L^1)$.  In this case we show that there is no
$\cW$--left approximate identity for $\cI$ (the same argument shows that it also has no $\cW$--right
approximate identity).

We use without further reference the facts  that $\mathcal{W} $ is equal to the ideal $\mathcal{S} =
\mathcal{S}(L^1)$ of strictly singular operators on $L^1$,  and that an operator on $L^1$ is weakly
compact if and only if it is $\ell^1$--singular; i.e., is not an isomorphism when restricted to any
subspace that is isomorphic to $\ell^1$. This background is presented in, for example, \cite[Chapter
5]{AK}.

Since $L^1$ has the Dunford--Pettis property \cite[Theorem 5.4.6]{AK}, the product $ \mathcal{W}
\cdot \mathcal{I }$ is a subset of the compact operators $\mathcal{K} = \mathcal{K}(L^1)$. Thus if
$(S_\af)$ is a net of weakly compact operators and $\| S_\af T - T\| \to 0$, then $T$ must be a
compact operator.  This proves, in particular that
there is no $\cW$--left approximate identity for $\cI$ unless
 $\cI \subset \cK$. This finishes the proof because $\cK$ is the smallest non trivial closed ideal
 in $L(L^1)$.

 Next suppose that $\cI$ is the ideal $\cW$ of weakly compact operators on $L^1$.  We have already
 observed that $\cW$ does not have a left approximate identity (since there are weakly compact
 operators on $L^1$ that are not compact), but we want to show that there is no $\cM$--left
 approximate identity for $\cW$, where
$\cM=\cM(L^1)$ is the unique maximal ideal in $L(L^1)$ described in the introduction. It then
follows that there is no $\cM$--left approximate identity for any closed ideal in $L(L^1)$ that
contains $\cW$. We then complete the proof of the $p=1$ case of Theorem \ref{main} by recalling
that a closed ideal  in $L(L^1)$ that is not contained in $\cW$ must contain $\cW$.  This is true
because (1) the identity operator $I_{\ell^1}$ on $\ell^1$  factors through every non weakly compact
operator in $L(L^1)$ \cite{KP}, and (2) every weakly compact (even every representable) operator
from $L^1$ to any Banach space factors through $\ell^1$. Again, this background is in the above
mention chapter in \cite{AK}.

To prove that  there is no $\cM$--left approximate identity for $\cW$,
 it is enough to find $S\in \cW$ such that  $S$ is not in the closure of $\cM \cdot S$.  To do this,
 first observe that $S:=I_{2,1}Q$ is in $\cW$, where $I_{2,1}$ is the formal identity from $L^2$
 into $L^1$ and $Q$ is a quotient mapping from $L^1$ onto $L^2$; that is, $Q$ maps the open unit
 ball $B^\circ_{L^1}$ of $L^1$ onto the open unit ball $B^\circ_{L^2}$ of $L^2$.  Thus
 $\overline{SB^\circ_{L^1}}$ contains all measurable functions having constant modulus $1$. But if
 $T$ is an operator on $L^1$ that  maps all measurable functions having constant modulus $1$
 into functions whose norms are bounded away from zero, then $T$ is not in $\cM$ by \cite{es},
 \cite{ros2}. In particular, no operator in $\cM$ can be close to the identity on the functions that
 have constant modulus $1$. This shows that $S$ is not in the closure of $\cM \cdot S$.

 This completes the proof of Theorem \ref{main} in the case $p=1$. We just remark that we do not
 really need the fact that a closed ideal $\cI$  in $L(L^1)$ that contains a non weakly compact
 operator must contain $\cW$, which uses the background facts (1) and (2).  Only (1) is needed. (1)
 says that  $\cI$ contains a projection whose range is isomorphic to $\ell^1$, and it is  easy to
 construct an operator $U: \ell^1 \to L^1$ such that $UB_{\ell_1}$ contains all measurable functions
 that have constant modulus one.

\bigskip\noindent
Case 2: $2<p<\infty$.

\bigskip
We assume that $\cI$ is a non trivial ideal in $L(L^p)$ different from $\cK=\cK(L^p)$.
The main tool needed for the proof is the Kadec--Pe\l czy\'nski dichotomy principle \cite{KP}, which
says that every normalized weakly null sequence $\xn$ in $L^p$, $2<p<\infty$, contains a subsequence
$\xn$ that is equivalent either to the unit vector basis of $\ell^2$ or to the unit vector basis of
$\ell^p$, and moreover that the closed linear space of $\xn$ is complemented in $L^p$. (This and
other standard background on $L^p$ that we use  is contained in the book \cite{AK}.) Both
possibilities can occur because $L^p$ contains complemented subspaces isomorphic to $\ell^2$ and
contractively complemented subspaces isometrically isomorphic to $\ell_p$.  By the K--P dichotomy,
if $T$ is a non compact operator on $L^p$ then either $i_p :=i_{p,p}$, $i_2:=i_{2,2}$, or $i_{2,p}$
factors through $T$, where for $r\le s$ we denote by $i_{r,s}$ the formal identity operator from
$\ell^r$ into $\ell^s$. If $T$ lies in the ideal $\cI$ and $i_{2,p}$ does not factor through $T$,
then the ideal property implies that $i_{2,p}$ factors through some other operator in $\cI$. So in
all cases we have $T\in \cI$ such that $i_{2,p}$ factors through
 $T$. However, we get more information when either $i_2$ or $i_p$ factors through an operator in
 $\cI$, so we discuss all three cases.

 The easiest case is that $i_2$ factors through some $T\in \cI$.  Thus in $TB_{L^p}$ there is a
 sequence $\xn$ that is equivalent to the unit vector basis for $\ell^2$ and such that the closed
 linear span $[\xn]$ of $\xn$ is complemented in $L^p$.  Let $\hn$ be the Haar basis for $L^p$,
 normalized so that each $h_n$ has norm one.  Since $\hn$ is an unconditional basis for $L^p$  and
 $L^p$ has type $2$, the operator $x_n\mapsto h_n$ extends to a bounded linear operator from $[\xn]$
 into $L^p$.  Since $[\xn]$ is complemented in $L^p$, this operator extends to an operator $S$ on
 $L^p$. Thus the operator $U:=ST$ is an operator in $\cI$ such that $UB_{L^p}$ contains the Haar
 basis $\hn$. We claim that $U$ is not in the closure of $\cM\cdot U$ (where here $\cM=\cM(L^p)$).  This claim is an immediate
 consequence of Andrew's result \cite{and}   that if $S\in L(L^p)$ and $\sup_n\|Sh_n - h_n\| <1$,
 then $S$ is not in $ \cM$. This proves that there is no $\cM$--left approximate identity for $\cI$.

Next assume that $i_p$ factors through some $T\in\cI$.  It then follows from Proposition
\ref{prp:newoperator} below that there is another operator $U\in \cI$ such that   $UB_{L^p}$
contains the Haar basis $\hn$. The argument in the previous paragraph then gives that there is no
$\cM$--left approximate identity for $\cI$.

 Finally, assume that only $i_{2,p}$ factors through an operator $T\in \cI$.  In this case we only
 prove that if $\cJ$ is a closed ideal such that there is a $\cJ$--left approximate identity for
 $\cI$, then $i_p$ factors through some operator in $\cJ$.  We can assume that there is $\xn$ in
 $TB_{L^p}$ so that $\xn$ is equivalent to the unit vector basis of $\ell^p$ and there is a
 projection $P$ from $L^p$ onto $[\xn]$.       Now suppose that there is a $\cJ$--left approximate
 identity for $\cI$.  Given $\ep>0$, there is $S\in \cJ$ such that $\sup_n \|x_n - Sx_n\| < \ep$.
 Replacing $S$ with $PS\in \cJ$, we can (since $\ep>0$ is arbitrary) assume that $S[\xn] \subset
 [\xn]$ and so the restriction of $S$ to $[\xn]$ is equivalent to a necessarily non compact operator
 on $\ell_p$ and hence $i_p$ factors through $S$.

 This completes the proof of the case $2<p<\infty$.  Notice that by the K-P result, an operator on
 $L^p$ is strictly singular if and only if neither $i_p$ nor  $i_2$ factors through the operator.
 Thus any closed ideal for which there is a $\cM$--approximate identity must be contained in the
 ideal of strictly singular operators.  However, we do not know whether any of these
 $2^{2^{\aleph_0}}$ closed  ideals other than $\cK$ has a $\cM$--left approximate identity.

\bigskip
Case 3: $1<p<2$.

\bigskip
Let $\cI$ be any non trivial ideal in $L(L^p)$ other than $\cK$.  We show that $\cI$ does not have a
$\cM$--left approximate identity. By duality from the $p>2$ case we know that there is $T\in \cI$
such that $i_{p,2}$ factors through $T$. Let $\wn$ be the Walsh functions. The mapping
$\delta_n\mapsto w_n$ extends to a contraction  from $\ell^2$ into $L^p$ (here $(\deltan)$ is the
unit vector basis for $\ell^2$).  Thus there is another operator $S\in \cI$ such that
$\wn \subset SB_{L^p}$.  Thus if there is a $\cM$--left approximate identity for $\cI$, then for
each $\ep>0$ there is $U\in \cM$ such that $\sup_n \| w_n - Uw_n\| < \ep$.  Another result from
Andrew's paper \cite{and} says that such a $U$ cannot be $L^p-$-singular if $\ep <1$.  But it was
proven in \cite[Theorem 9.1]{jmst} that no non $L^p$--singular operator on $L^p$ is in $\cM$ (this
just says that an operator on $L^p$ that preserves some isomorph of $L^p$ must send some isomorph of
$L^p$ isomorphically onto a complemented subspace of $L^p$). This completes the proof of the case
$1<p<2$, but
we remark that one could avoid quoting the result from \cite{jmst} by understanding what Andrews
proved but did not state explicitly in \cite{and}.
\qed

By using duality we can summarize the main things that were proved in Theorem \ref{main} as
Corollary \ref{cor}. In addition to background that has already been stated, we need the fact
\cite{Weis} that an operator in $L(L^p)$, $1<p<2$,     is $L^p$--singular if it is both
$\ell^p$--singular and $\ell^2$--singular because this fact does not follow formally by duality from
the corresponding statement for $2<p<\infty$.

\begin{cor}\label{cor}
Let $\cI$ be a nontrivial closed ideal in $L(L^p)$, $1\le p <\infty$, other than the ideal of
compact operators. Then $\cI$ does not have a left approximate identity. When $p>1$, $\cI$ does not
have a right approximate identity. Moreover, if $1<p<2$, then $\cI$ does not have a $\cM$--left
approximate identity and when $2<p<\infty$,  $\cI$ does not have a $\cM$--right approximate
identity. If $1<p<\infty$ and $\cI $ is not contained in the ideal of strictly singular operators,
then $\cI$ does not have either a $\cM$--left approximate identity or a $\cM$--right approximate
identity.
If $p=1$ and $\cI $ is not contained in the ideal of weakly compact operators, then $\cI$ does not
have a either a $\cM$--left approximate identity or a $\cM$--right approximate identity.
\end{cor}

\bigskip
In the proof of Theorem \ref{main} we used the fact that there exists $T\in L(\ell^p,L^p)$,
$2<p<\infty$, such that $TB_{\ell^p}$ contains the $L^p$ normalized Haar basis for $L^p$.  This fact
is immediate from Proposition \ref{prp:newoperator} and the well known fact that for $1<p<\infty$
the space $\ell^p$ is isomorphic to $(\sum_{n=1}^\infty \ell^2_n)_p$, or to any other  $p$--sum of a
sequence of  finite dimensional non zero Hilbert spaces  \cite{KP}.

Recall the definition of the Haar functions: For $n=0,1,\dots$ and $i=1,\dots,2^n$, denote
\[
h_{n,i}(t)=\begin{cases}
1& \frac{i-1}{2^n}\le t < \frac{2i-1}{2^{n+1}}\\
-1&\frac{2i-1}{2^{n+1}}\le t < \frac{i}{2^{n}}.
\end{cases}
\]
As is well known, the Haar system, ${\bf 1}\cup\{h_{n,i}\}_{n=0,i=1}^{\infty\ \ \ 2^n}$, in its
lexicographic order, forms a Schauder basis for $L^p$, $1\le p<\infty$ and this basis is
unconditional iff $p>1$.

We also introduce the space $(\sum_{n=0}^\infty\oplus\ell_{2^n}^2)_p$; i.e., the $\ell_p$ sum of
Euclidean spaces of dimensions $2^n$, $n=0,1,\dots$. Its natural basis will be denoted
$\{e_{n,i}\}_{n=0,i=1}^{\infty\ \ \ 2^n}$. In the next proposition we show that, for $p>2$, the
basis to basis map $e_{n,i}\to h_{n,i}/\|h_{n,i}\|_p$ extends to a bounded operator.

\begin{prp}\label{prp:newoperator}
Let $2<p<\infty$ and define $T$ from the finitely supported sequences in
$(\sum\oplus\ell_{2^n}^2)_p$ to $L^p$ as the linear extension of $Te_{n,i}=h_{n,i}/\|h_{n,i}\|_p$.
Then $T$ is bounded and thus naturally extends to an operator from $(\sum\oplus\ell_{2^n}^2)_p$ to
$L^p$.
\end{prp}

\pf
Recall (\cite{mu}, 1.2.6, page 37) that dyadic BMO  is the space of all formal expansions
$\sum_{n,i}b_{n,i}h_{n,i}$ whose norm
\[
\|\sum_{n,i}b_{n,i}h_{n,i}\|_{BMO}=\sup_{0\le n<\infty,1\le i\le 2^n}(2^n\sum_{m=n}^\infty\sum_{\{j;
|h_{m,j}|\le |h_{n,i}|\}}b_{m,j}^22^{-m})^{1/2}
\]
Is finite. So,
\begin{multline*}
\|\sum_{n,i}b_{n,i}h_{n,i}\|_{BMO}\le \sup_{0\le n<\infty}(2^n\sum_{m=n}^\infty
2^{-m}\sum_{j=1}^{2^m}b_{m,j}^2)^{1/2}\\
\le
\sup_{0\le n<\infty}(2\sup_{m\ge
n}\sum_{j=1}^{2^m}b_{m,j}^2)^{1/2}=2^{1/2}\|\{b_{m,j}\}\|_{(\sum\oplus\ell_{2^n}^2)_\infty}.
\end{multline*}
Recall that, for $1\le p<\infty$, dyadic $H^p$ is the space of all functions
$\sum_{n,i}a_{n,i}h_{n,i}$ for which
\[
\|\sum_{n,i}a_{n,i}h_{n,i}\|_{H^p}=\|(\sum_{n,i}a_{n,i}^2h_{n,i}^2)^{1/2}\|_{H^p}
\]
is finite.

By the duality between dyadic $H^1$ and dyadic BMO (see section 1.2 in \cite{mu}) we get that for
some universal constant $C$ and all coefficients $\{a_{n,i}\}$,
\begin{multline}\label{eq:h1}
\sum_{n=0}^\infty(\sum_{i=1}^{2^n}a_{n,i}^2)^{1/2}\\
\le
C \|\sum_{n,i}a_{n,i}h_{n,i}/\|h_{n,i}\|_{L^1}\|_{H^1}
=
 C\|(\sum_{n,i}a_{n,i}^2h_{n,i}^2/\|h_{n,i}\|_{L^1}^2)^{1/2}\|_{L^1}.
\end{multline}
We now interpolate between $H^1$ and $L^2$ to get that for some absolute constant $C$, for all
$1<p<2$ and all coefficients $\{a_{n,i}\}$,
\begin{multline}\label{eq:hp}
\left(\sum_{n=0}^\infty(\sum_{i=1}^{2^n}a_{n,i}^2)^{p/2}\right)^{1/p}\\
\le
C\left\|\sum_{n,i}a_{n,i}h_{n,i}/\|h_{n,i}\|_{L^p}\right\|_{H^p}
=
C \left\|\big(\sum_{n,i}a_{n,i}^2h_{n,i}^2/\|h_{n,i}\|_{L^p}^2\big)^{1/2}\right\|_{L^p}.
\end{multline}
More precisely, we use
\begin{Thm*}[2.3.4 in \cite{mu}] For all $v=\sum_{n,i}v_{n,i}h_{n,i}\in H^p, 1<p<2$, there exists
$x=\sum_{n,i}x_{n,i}h_{n,i}\in H^1$ \ and \ $y=\sum_{n,i}y_{n,i}h_{n,i}\in L^2$ \ satisfying
$x_{n,i},y_{n,i}\ge 0$ for all $n,i$,
\[
|v_{n,i}|=x_{n,i}^{1-\theta}y_{n,i}^\theta \ \ \mbox{for all} \ \ n=0,1,\dots, \ \mbox{and}\
i=1,\dots,2^n
\]
and
\[
\|x\|_{H^1}^{1-\theta}\|y\|_{L^2}^\theta\le 32\|v\|_{H^p},
\]
where $\theta$ is defined by $\frac{1}{p}=\frac{1-\theta}{1}+\frac{\theta}{2}$; i.e.,
$\theta=2-\frac{2}{p}$.
\end{Thm*}
 To prove (\ref{eq:hp}) we use this theorem with $v_{n,i}=a_{n,i}/\|h_{n,i}\|_{L^p}=2^{n/p}a_{n,i}$.
 So there are non negative $x_{n,i},y_{n,i}$ such that
 \[
 2^{n/p}|a_{n,i}|=x_{n,i}^{1-\theta}y_{n,i}^\theta \ \ \mbox{for all} \ \ n=0,1,\dots, \ \mbox{and}\
 i=1,\dots,2^n
\]
and
\begin{equation}\label{eq:int}
\|(\sum_{n,i}x_{n,i}^2 h_{n,i}^2)^{1/2}\|_{L^1}^{1-\theta}\|(\sum_{n,i}y_{n,i}^2
h_{n,i}^2)^{1/2}\|_{L^2}^{\theta}\le 32\|\sum_{n,i}a_{n,i}h_{n,i}/\|h_{n,i}\|_{L^p}\|_{H^p}.
\end{equation}
Now,
\begin{align*}
&(\sum_{n=1}^\infty(\sum_{i=1}^{2^n}a_{n,i}^2)^{p/2})^{1/p}=(\sum_n 2^{-n}(\sum_i
x_{n,i}^{2(1-\theta)}y_{n,i}^{2\theta})^{p/2})^{1/p}\\
&\phantom{aaaaa}\le
(\sum_n 2^{-n}(\sum_i x_{n,i}^2)^{\frac{(1-\theta)p}{2}}(\sum_i y_{n,i}^2)^{\frac{\theta
p}{2}})^{1/p}\ \ \mbox{by H\"older with}\ \ \frac{1}{1-\theta},\frac{1}{\theta}\\
&\phantom{aaaaa}
=
(\sum_n (2^{-2n}\sum_i x_{n,i}^2)^{\frac{(1-\theta)p}{2}}(2^{-n}\sum_i y_{n,i}^2)^{\frac{\theta
p}{2}})^{1/p}\\
&\phantom{aaaaa}\le
(\sum_n (2^{-2n}\sum_i x_{n,i}^2)^{1/2})^{1-\theta}(\sum_n 2^{-n}\sum_i y_{n,i}^2)^{\frac{\theta
}{2}}\\
&\phantom{aaaaaaaaaaaaaaaaaaaaaaaaaaaaaaaaaaaaa}\ \ \mbox{by H\"older with}\ \
\frac{1}{(1-\theta)p},\frac{2}{\theta p}\\
&\le C^{1-\theta}
\|\sum_{n,i} 2^{-n}x_{n,i}h_{n,i}/\|h_{n,i}\|_{L^1}\|_{H^1}^{1-\theta}
\|\sum_{n,i} 2^{-n/2}y_{n,i}h_{n,i}/\|h_{n,i}\|_{L^2}\|_{L^2}^{\theta}\ \ \mbox{by}\ (\ref{eq:h1})\\
&= C^{1-\theta}
\|\sum_{n,i} x_{n,i}h_{n,i}\|_{H^1}^{1-\theta}
\|\sum_{n,i} y_{n,i}h_{n,i}\|_{L^2}^{\theta}\\
&<
32C \|\sum_{n,i}a_{n,i}h_{n,i}/\|h_{n,i}\|_{L^p}\|_{H^p}\ \ \mbox{by}\ (\ref{eq:int}).
\end{align*}
The unconditionality of the Haar system in $L^p$ implies that the last quantity is equivalent, with
constant depending on $p$ only, to
\[
\|\sum_{n,i}a_{n,i}h_{n,i}/\|h_{n,i}\|_{L^p}\|_{L^p}.
\]
Dualizing again we get the required inequality showing that $T$ is bounded.
\qed

\bigskip
Although we focus here on separable $L^p(\mu)$ spaces, it is worth remarking that if $L^p(\mu)$ is
non separable, then there is a non trivial ideal $\cI$ in $L(L^p(\mu))$  with $\cI \not=
\cK(L^p(\mu))$
such that $\cM$ has a left approximate identity consisting of contractive projections; namely, let
$\cI$ be the ideal of operators that have separable ranges.  Since every closed sublattice of $L^p$
is contractively complemented \cite[Lemma 1.b.9]{LT2}, the collection of contractive projections
onto separable sublattices forms a left approximate identity for $\cI$.

\section {Right approximate identities for ideals in $L(L^1)$}\label{L^1}

It was stated in Corollary \ref{cor} that no small ideal (that is, contained in the ideal of
strictly singular operators)  in $L(L^1)$) other than the compact operators has a $\cW (L^1)$-right
approximate identity. In this section we consider two of the three known large closed  ideals in
$L(L^1)$. We do not know whether the remaining large ideal, the ideal of Dunford-Pettis operators,
has a right approximate identity. Our guess is that it does not.

Background we use is contained in \cite[section III]{du}.  The main facts needed are that RNP
(representable) operators with domain an $L^1$ space are the same as the operators that factor
through $\ell^1$, that every operator $B$ from $L^1$ to $\ell^1$ is an RNP operator, and if $B$ is
an RNP operator with domain $L^1$, then there is a partition $(E_n)_{n=1}^\infty$ of $(0,1)$ such
that for each $n$ the operator $B_n := B_{|L^1(E_n)}$ is compact (this is how Lewis and Stegall
\cite{ls} showed that RNP operators from $L^1$ factor through $\ell^1$--that operators from $L^1$
into $\ell^1$ are RNP operators was known before their work).

The two large ideals in $L(L^1)$ that we consider are
 $\Gamma_{\ell^1}(L^1)$, the ideal of operators in $L(L^1)$ that factor through $\ell^1$; and
 $\cM(L^1)$, the unique maximal ideal in $L(L^1)$. The ideal $\cM(L^1)$ is the class of the
 operators $T$  in $L(L^1)$ such that the identity $I_{L^1}$ does not factor through $T$. It was
 already mentioned in the introduction that  $\cM(L^1)$ is indeed an ideal (so that it obviously is
 maximal and the only maximal ideal) and that $\cM(L^1)$ can also be described as the $L^1$-singular
 operators on $L^1$.

\begin{thm}\label{gamma1}
The ideal $\cI:=\Gamma_{\ell^1}(L^1)$ in $L(L^1)$ has a right contractive approximate identity
consisting of contractive idempotents.
\end{thm}

\pf Notice that $\cI$  is a large ideal and is closed because the RNP operators are closed.  To
prove that $\mathcal{I}$  has a right approximate identity consisting of contractive idempotents it
is enough to show that there is a net $\mathcal{N}$ of contractive idempotents (i.e., conditional expectations)
in $\mathcal{I}$  such that for all $B:L^1\to \ell^1$  and each $\epsilon >0$ there is a contractive
idempotent $P$ in $\mathcal{N}$ such that $\|BP - B\| \le \epsilon$.

Let $\mathcal{N}$ be the net of all $P\in L(L^1)$ such that there exist a partition $\{E_n\}_{n=1}^\infty$ of $(0,1)$
and, for each $n$, a finite rank idempotent  $P_n$ on $L^1(E_n)$ such that, for each $x\in L^1$ and each $n$, $(Px)_{|E_n}=P_n(1_{E_n}x)$.

Since each $P_n$, being of finite rank, is in $\mathcal{I}$ (formally, $P_n$ preceded by the restriction to $E_n$ and followed by the natural embedding of $L^1(E_n)$ into $L^1$) it is easy to see that $P$ is in $\mathcal{I}$ (for if $A_1,A_2,\dots $ all
factor through $\ell^1$, then so does $A:= (A_1,A_2,\dots, ) : L^1 \to (\sum_{k=1}^\infty L^1)_1$, where
$Ax :=(A_1x,A_2x,\dots)$). Also, $\mathcal{I}$ is clearly a net under the natural order on idempotents.

Finally, given a $B:L^1\to \ell^1$  and  $\epsilon >0$
take a partition $(E_n)$ as in the result of Lewis and Stegall mentioned above and define $B_n := B_{|L^1(E_n)}$.  Since $B_n$ is compact there
is a contractive finite rank idempotent $P_n$ on ${L^1(E_n)}$ such that
$\|B_n - B_nP_n\| \le \epsilon$. Define $P$ on $L^1$ by $(Px)_{|E_n} := P_n (1_{E_n}x)$.  
\qed

\medskip

Y. Choi \cite{choi} has used Theorem \ref{gamma1} to give a simple proof that $L(L^1)$ is not
amenable as a Banach algebra.

In order to see that $\cM(L^1)$ does not have a right approximate identity, we need some further
background and notation.
Let $r_1,r_2,\dots$, denote the Rademacher functions, considered either as functions on $[0,1]$ or
the coordinate projections on the Cantor group $\{-1,1\}^N$ endowed with Haar measure. The
($L^1$-normalized) Haar system (less the constant function), used already in the previous section,
can also be defined by:
\[
h_\phi=r_1,\quad\quad
h_{\d}=r_2(1+\d r_1),\quad \d=\pm1,
\]
and for $n=2,3,\dots$, by induction,
\[
h_{\d}=r_{n+1}\prod_{i=1}^n(1+\d_ir_i),\quad \d=(\d_1,\dots,\d_n),\quad \d_i=\pm1.
\]
This system union a non zero constant function spans $L^1$ and is a martingale difference sequence
in any order such that $(\d_1,\dots,\d_n)$ comes after $(\gamma_1,\dots,\gamma_{n-1})$ for every
$n>1$ and every $\d_1,\dots,\d_n$ and $\gamma_1,\dots,\gamma_{n-1}$ in $\{-1,1\}$ (here and
elsewhere when $n=1$,  $(\gamma_1,\dots,\gamma_{n-1})$ is interpreted as $\phi$). Hence the Haar
system forms a monotone basis of $L^1$ in any such order.

\begin{claim}\label{claim:haar} Let
$\varphi:\{\d_1,\dots,\d_n\}_{n=1,\d_i=\pm1}^\infty\to\{2,3,\dots\}$ be a function such that
\[
\varphi(\d_1,\dots,\d_n)>\varphi(\gamma_1,\dots,\gamma_{n-1})
\]
for every  $n>1$ and every $\d_1,\dots,\d_n$ and $\gamma_1,\dots,\gamma_{n-1}$ in $\{-1,1\}$. Define
\[
k_\phi=r_1,\quad\quad
k_{\d}=r_{\varphi(\d)}(1+\d r_1),\quad \d=\pm1,
\]
and for $n=2,3,\dots$, by recursion,
\[
k_{\d}=r_{\varphi(\d_1,\dots,\d_n)}\prod_{i=1}^n(1+\d_ir_{\varphi(\d_1,\dots,\d_{i-1})}),\quad
\d=(\d_1,\dots,\d_n),\quad \d_i=\pm1.
\]
Then   $k_\phi,\{k_{(\d_1,\dots,\d_n)}\}$ has the same joint distribution as the Haar system
$h_\phi,\{h_{(\d_1,\dots,\d_n)}\}$.
\end{claim}

\pf Define a map from the $\sigma$-field generated by $r_1,\dots,r_{n+1}$ into (but not necessarily
onto) the $\sigma$-field generated by $\{r_1,r_{\varphi(\eta_1)},\dots,r_{\varphi(\eta_1,\dots,\eta_n)}\}_{\eta_1,\dots,\eta_n\in\{-1,1\}}$ by
\begin{multline*}
\psi(\{t;\ r_1(t)=\d_1,\dots,r_n(t)=\d_n,r_{n+1}(t)=\d_{n+1}\})\\
=\{t;\ r_1(t)
=\d_1,r_{\varphi(\d_1)}=\d_2,\dots,r_{\varphi(\d_1,\dots,\d_n)}(t)=\d_{n+1}\}.
\end{multline*}
This is a measure preserving map and $k_{(\d_1,\dots,\d_m)}(\psi(t))=h_{(\d_1,\dots,\d_m)}(t)$ for
all $t$, $m=1,\dots,n$ and $\d_1,\dots,\d_n\in\{-1,1\}$.
 \qed

\begin{lem}\label{lemma:rai}
Let $\{x_k\}_{k=1}^\infty$ be a sequence in $L^1$ such that there is a sequence
$\{x^*_k\}_{k=1}^\infty$ of norm one functionals  with $x^*_k(x_k)\ge a>0$ for all $k$ and
$x^*_k(x_l)=0$ for $k\not=l$.
Let $\{h_k\}_{k=1}^\infty$ be an enumeration of the Haar system in $L^1$ and assume that
$S(h_k)=x_k$ extends to an   operator in $L(L^1)$. Assume also that $T$ in $L(L^1)$ satisfies
\[
\|ST-S\|=b<a.
\]
Then $T$ is an Enflo operator; that is, there is a subspace $Y$ of $L^1$ so that the restriction
$T_{|Y}$ of $T$ to $Y$ is an  isomorphism.
\end{lem}
\pf
By the assumptions,
\[
x^*_k(STh_k)\ge x^*_k(x_k)-b\ge a-b.
\]
Write
\[
Th_k=\alpha h_k+\sum_{l\not= k}
\alpha_l h_l
\]
so that
\[
x^*_k(STh_k)=\alpha x^*_k(Sh_k)\le |\alpha|\|S\|
\]
and hence
\[
|\alpha|\ge (a-b)/\|S\|.
\]
It now follows from \cite{lmms} that $T$ is an Enflo operator.
\qed

For $0<\e<1$,   the biased coin operator $T_\e$  is convolution with the measure
$\prod_{i=1}^\infty(1+\e r_i)$. The basic source for information about biased coin operators is
Bonami's thesis \cite{Bon}.  Biased coin operators  were introduced  into Banach space theory by
Rosenthal \cite{Ros1}.

\begin{thm}\label{theorem:norai}
Suppose that $\cI$ is a closed ideal in $L(L^1)$ and there is $0<\e<1 $ such that the biased coin
operator $T_\e $ is in $\cI$. Then $\cI$ does not have a  $\cM(L^1)$--right approximate identity. In
particular, the largest ideal $\cM(L^1)$ does not have a right approximate identity.
\end{thm}

\pf Recall that $\cM(L^1)$ is the collection of non Enflo operators. $T_\e$ is not an Enflo operator
\cite{Ros1}, so $T_\e$ is contained in  $\cM(L^1)$; so the ``In particular" statement is clear. To
prove the Theorem itself,
by Lemma \ref{lemma:rai} it is enough to find an operator $S$ in  $\cI$ that sends the Haar system
to a biorthogonal sequence whose biorthogonal functionals are bounded. Note that $T_\e
h_{(\d_1,\dots,\d_n)}=\e r_{n+1}\prod_{n=1}^n(1+\e \d_i r_i)$, and $\|T_\e
h_{(\d_1,\dots,\d_n)}\|=\e$ is indeed bounded away from zero. However, this sequence does not have a
bounded sequence of biorthogonal functionals. We overcome this difficulty by considering the action
of $T_\e$ on an equivalent system of the kind discussed in Claim \ref{claim:haar}.

Let the function $\varphi:\{\d_1\dots,\d_n\}_{n=1,\d_i=\pm1}^\infty\to\{2,3,\dots\}$  in the
statement of Claim \ref{claim:haar}  also be one to one. Note that
\[
T_\e k_{\d}=\e r_{\varphi(\d_1,\dots,\d_n)}\prod_{i=1}^n(1+\e \d_ir_{\varphi(\d_1,\dots,\d_{i-1})}).
\]
It is still the case that $\|T_\e k_{\d}\|=\e$,  but now, due to the fact that
$\varphi(\d_1,\dots,\d_n)$ are all different, the sequence $\{T_\e k_{\d}\}$ is a martingale
difference sequence and thus a montone basic sequence.  Consequently, the biorthogonal functionals
to $\{T_\e k_{\d}\}$ all have norm at most $2/\e$.  Thus we can apply   Lemma \ref{lemma:rai}, using
for
 the operator $S$  the basis to basis map from the Haar system to the $k$-system followed by $T_\e$.
\qed

\noindent William B. Johnson\newline
             Department Mathematics\newline
             Texas A\&M University\newline
             College Station, TX, USA\newline
             E-mail: johnson@math.tamu.edu

\bigskip
\noindent Gideon Schechtman\newline Department of
Mathematics\newline Weizmann Institute of Science\newline
Rehovot, Israel\newline E-mail: gideon.schechtman@weizmann.ac.il

\end{document}